
\magnification 1200
\parskip 9pt plus 5pt minus 3pt
\def\item{\parskip 5pt plus 3pt minus 3pt \par\hang\textindent}

\headline{\ifodd\pageno\ifnum\pageno=1\else\hfil\rlap{\hautdepage\folio}\fi\else\llap{\hautdepage\folio}\hfil\fi} 
\footline{}              
\catcode`\@=11
\newdimen\margereliured \margereliured=1truecm
\newdimen\margereliureg \margereliureg=-0.1truecm
\def\plainoutput{%
\shipout\vbox{
\ifodd\pageno \moveright\margereliured
              \else \moveleft\margereliureg \fi
\hbox{\vbox{\makeheadline\pagebody\makefootline}}}%
\advancepageno
\ifnum\outputpenalty>-\@MM \else\dosupereject\fi}
\catcode`\@=12

\def\tenpoint{%
  \textfont0=\tenrm \scriptfont0=\sevenrm \scriptscriptfont0=\fiverm
  \def\rm{\fam0\tenrm}%
  \textfont1=\teni \scriptfont1=\seveni \scriptscriptfont1=\fivei
  \def\oldstyle{\fam1\teni}%
  \textfont2=\tensy \scriptfont2=\sevensy \scriptscriptfont2=\fivesy
  \textfont\itfam=\tenit
  \def\it{\fam\itfam\tenit}%
  \def\sl{\fam\slfam\tensl}%
  \textfont\slfam=\tensl \scriptfont\slfam=\sevensl \scriptscriptfont\slfam=\fivesl
  \def\bf{\fam\bffam\tenbf}%
  \textfont\bffam=\tenbf \scriptfont\bffam=\sevenbf
  \scriptscriptfont\bffam=\fivebf  
  \def\tt{\fam\ttfam\tentt}%
  \textfont\ttfam=\tentt
  \abovedisplayskip=6pt plus 2pt minus 6pt
  \abovedisplayshortskip=0pt plus 3pt
  \belowdisplayskip=6pt plus 2pt minus 6pt
  \belowdisplayshortskip=7pt plus 3pt minus 4pt
  \smallskipamount=3pt plus 1pt minus 1pt
  \medskipamount=6pt plus 2pt minus 2pt
  \bigskipamount=12pt plus 4pt minus 4pt
  \normalbaselineskip=12pt
  \setbox\strutbox=\hbox{\vrule height8.5pt depth3.5pt width0pt}%
  \normalbaselines\rm}
\catcode`\|=13
\def\today{\ifcase\month\or january \or february \or march \or april
\or may \or june\or july\or august \or september\or october\or november\or
December\fi\ \number\day , \number\year}


\newskip\LastSkip
\def\nobreakatskip{\relax\ifhmode\ifdim\lastskip>0pt
  \LastSkip\lastskip\unskip
  \nobreak\hskip\LastSkip
  \fi\fi}
\catcode`\;=\active \def;{\nobreakatskip\string;}
\catcode`\:=\active \def:{\nobreakatskip\string:}
\catcode`\!=\active \def!{\nobreakatskip\string!}
\catcode`\?=\active \def?{\nobreakatskip\string?}

\newif\ifrefvis
\refvisfalse

\newif\ifarxiv
\arxivfalse
\def\arXiv{\arxivtrue}

\newif\iffrance

\def\anglais{\francefalse}

\newskip\afterskip
\catcode`\@=11
\def\p@int{.\par\vskip\afterskip\penalty100} 
\def\p@intir{\discretionary{.}{}{.\kern.35em---\kern.7em}}
\def\pointir{\afterassignment\pointir@\global\let\next=}
\def\pointir@{\ifx\next\par\p@int\else\p@intir\fi\egroup\next}
\catcode`\@=12
\def|{\relax\ifmmode\vert\else\findef\fi}
\def\findef{\errhelp{Cette barre verticale ne correspond ni a un \vert mathematique
                        ni a une fin de definition, le contexte doit vous indiquer ce qui manque.
                        Si vous vouliez inserer un long tiret, le codage recommande est ---,
                        dans tous les cas, la barre fautive a ete supprimee.}%
                        \errmessage{Une barre verticale a ete trouvee en mode texte}}

\def\TITR#1|{\null{\mss\baselineskip=17pt
                           \vskip 3.25ex plus 1ex minus .2ex
                           \leftskip=0pt plus \hsize
                           \rightskip=\leftskip
                           \parfillskip=0pt
                           \noindent #1
                           \par\vskip 2.3ex plus .2ex}}
 
\def\auteur#1|{\penalty 500
               \vbox{\centerline{
                 \iffrance par \else by \fi #1}
                \vskip 10pt}\penalty 500}


\def\fonction#1|#2|#3|#4|{\lower 6pt\hbox{
$\matrix{#1 \hfill &\longrightarrow \hfill & #2\hfill\cr
\hfill #3  &\longmapsto \hfill &#4 \hfill\cr}$}}
\newcount\thesection
\newcount\thesubsection
\newcount\thesubsubsection
\newcount\theparagraf
\newcount\thetheo
\newcount\theequ
\global\thesection=0
\global\thesubsection=0
\global\thesubsubsection=0
\global\theparagraf=0
\global\thetheo=0
\global\theequ=0

\font\sevensl= cmsl7
\font\fivesl= cmsl5

\font\tentite = cmbx10 at 16pt
\font\seventite = cmbx7 at 11pt
\font\fivetite = cmbx5 at 8pt
\newfam\titefam
\textfont\titefam = \tentite
\scriptfont\titefam = \seventite
\scriptscriptfont\titefam = \fivetite

\font\sectfont = cmbx10 at 14pt
\font\sectscript = cmbx7 at 10pt
\font\sectsscript = cmbx5 at 7pt
\newfam\sectionfam 
\textfont\sectionfam = \sectfont 
\scriptfont\sectionfam = \sectscript
\scriptscriptfont\sectionfam = \sectsscript
\def\sectionfont{\fam\sectionfam\sectfont}

\font\subsectfont =  cmbx10 at 12pt
\font\subsectscript = cmbx7 at 8pt
\font\subsectsscript = cmbx5 at 6pt
\newfam\subsectionfam
\textfont\subsectionfam = \subsectfont
\scriptfont\subsectionfam = \subsectscript
\scriptscriptfont\subsectionfam = \subsectsscript
\def\subsectionfont{\fam\subsectionfam\subsectfont}

\font\mss=cmss12 scaled \magstep1

\font\hautdepage = cmss8
\font\sf =cmss10
\font\ding = pzdr at 10pt

\font\tenmsb=msbm10
\font\sevenmsb=msbm7
\font\fivemsb=msbm5
\newfam\msbfam
\textfont\msbfam=\tenmsb
\scriptfont\msbfam=\sevenmsb
\scriptscriptfont\msbfam=\fivemsb
\def\Bbb#1{{\fam\msbfam\relax#1}}

\font\teneufrak=eufm10
\font\seveneufrak=eufm7
\font\fiveeufrak=eufm5
\newfam\eufrak
\textfont\eufrak=\teneufrak
\scriptfont\eufrak=\seveneufrak
\scriptscriptfont\eufrak=\fiveeufrak
\def\mathfrak#1{{\fam\eufrak\relax#1}}


\def\begincentered{\par\begingroup
\def \par{\hss\egroup\line\bgroup\hss}\obeylines
\line\bgroup\hss}
\def\endcentered{\hss\egroup\endgroup}

\hsize=12.5cm
\vsize=20cm
\parindent=1cm
\baselineskip=13pt
\hoffset=-0.1cm
\voffset=0.5cm

\long\def\partie#1{\begingroup\sectionfont
        \par\penalty -500
        \vskip 3.25ex plus 1ex minus .2ex
        \skip\afterskip=1.5ex plus .2ex
        \baselineskip=17pt        
        \par
        \def \par{\hss\egroup\line\bgroup\hss}\obeylines
        \line\bgroup\hss
\global\advance\thesection by 1 
\xdef \lastref{\number\thesection}
\global\thesubsection=0 \global\theparagraf=0  \number\thesection\quad#1\hss\egroup\endgroup\par}

\long\def\subpartie#1{\begingroup\subsectionfont%
                          \par\penalty -200
                          \vskip 3.25ex plus 1ex minus .2ex
                          \skip\afterskip=1.5ex plus .2ex
                          \baselineskip=15pt
        \par
        \def \par{\hss\egroup\line\bgroup\hss}\obeylines
        \line\bgroup\hss
\global\advance\thesubsection by 1
\xdef \lastref{\number\thesection.\number\thesubsection}
\global\thesubsubsection=0 \global\theparagraf=0
\number\thesection.\number\thesubsection\quad #1\hss\egroup\endgroup \par}

\def\\{\par}

\long\def\subsubpartie#1{\bgroup\bf
                                \par\penalty -100
                                \vskip 3.25ex plus 1ex minus .2ex
                                \skip\afterskip=1.5ex plus .2ex
                                
\begincentered
\global\advance\thesubsubsection by 1
\xdef \lastref{\number\thesection.\number\thesubsection.\number\thesubsubsection}
\global\theparagraf=0
\number\thesection.\number\thesubsection.%
\number\thesubsubsection\quad #1\endcentered\egroup\par}
\def\paragrafsubsub{%
\global\advance\theparagraf by 1
\xdef \lastref{\number\thesection.\number\thesubsection.\number\thesubsubsection.%
\romannumeral\theparagraf}
{\bf \number\thesection.\number\thesubsection.\number\thesubsubsection.\romannumeral\theparagraf}%
\kern0.2em --- }

\def\paragrafsub{%
\global\advance\theparagraf by 1
\xdef \lastref{\number\thesection.\number\thesubsection.\number\theparagraf}
{\bf \number\thesection.\number\thesubsection.\number\theparagraf}%
\kern0.2em --- }

\def\paragraf{%
\global\advance\theparagraf by 1 
\xdef \lastref{\number\thesection.\number\theparagraf}
{\bf \number\thesection.\number\theparagraf}%
\kern0.2em --- }

\def\paragraphe{%
\par \indent
\ifcase\thesubsection %
  \paragraf
\else
\ifcase\thesubsubsection\paragrafsub %
 \else\paragrafsubsub\fi
\fi}

\def\numeq{\global\advance\theequ by 1%
        \xdef \lastref{(\number\theequ)}%
        \eqno{(\number\theequ)}}
\def\nume{\global\advance\theequ by 1 
        \xdef \lastref {(\number\theequ)}%
        (\number\theequ)}
\def\theo{\global\advance\thetheo by 1%
\xdef\lastref{\number\thetheo}%
\number\thetheo}

\newwrite\fileref
\newread\instream
\newif \ifrefmodif 
\def \defineref#1#2{{\def\next{#1}%
        \expandafter\xdef
           \csname ref = \meaning\next\endcsname{#2}%
        }}

\def \initfileref{
         \openin\instream=\jobname.ref 
         \ifeof\instream \message{ Ahem } 
                \message{**********************************************************}
                \message{******* Le fichier \jobname.ref n'existait pas ***********}
                \message{********** Il faudra absolument recompiler ***************}
                \message{**********************************************************}
                \message{ }
         \else \closein\instream \input \jobname.ref \refmodiffalse 
         \fi 
        \ifarxiv\else\immediate\openout \fileref=\jobname.ref\fi
        \global\let\initfileref=\relax} 

\def \label#1{\ifarxiv\initfileref\else{%
        \toks0={#1}\wlog{REF \the\toks0= \lastref}
        \initfileref
        \immediate\write\fileref{\noexpand\defineref%
                {\the\toks0}{\lastref}}%
        \def\next{#1}%
        \expandafter\ifx 
          \csname ref = \meaning\next\endcsname\lastref
        \else \global\refmodiftrue  \message{ }\message{ Attention }
        \message{reference {\the\toks0 = \lastref }  modifiee ou redefinie} \message{ }\fi 
        \defineref{#1}{\lastref}%
        \ifrefvis\ifhmode\raise 6pt \vbox{\hbox to 0pt{\hss\fivebf [\the\toks0]\hss}}
        \else\llap{\fivebf [\the\toks0]}\fi\fi}\fi}

\def \ref#1{{\initfileref \def\next{#1}%
        \csname ref = \meaning\next\endcsname}}

\outer\def \bye{
        \closeout \fileref \vfill \supereject
        \ifrefmodif \erreurmodif \fi
        \end}

\def \erreurmodif{
        \message{ }
        \message{**********************************************************}
        \message{BEWARE, some references have been modified.}%
        \message{Please compile the tex file again to get references right,}
        \message{if this message appear again, then a reference must have}
        \message{been define at least two times.}        
        \message{**********************************************************}
        \message{ }}


\def \ignorepar{\afterassignment\ignoreparaux \let\next=}
\def \ignoreparaux{\ifx\next\par \let\next\ignorepar \fi \next}

\def\bibliography{\bgroup\sectionfont
        \par\penalty -500
        \vskip 3.25ex plus 1ex minus .2ex
        \skip\afterskip=1.5ex plus .2ex
\begincentered References \endcentered\egroup\par
\tabskip5pt minus 1pt%
\noindent\halign to \hsize \bgroup##\hfil&##\hfil\cr}
\def\bibliend{\egroup}

\def\bibitem[#1]#2{\xdef\lastref{[#1]}\label{#2}\hskip-8pt[#1]&%
\vtop\bgroup\hsize10.6cm\noindent\ignorespaces}
\def\finitem{\medskip \egroup \cr}

\def\newblock{\hskip .11em plus .33em minus .07em}
\def\cite#1{\ref{#1}}

\def\em{\sl}

\def\proof{\noindent{\bf Proof. }}
\def\qed{\ \-\hfill {\ding \char 111} \par}


\def\T{{\Bbb T}}
\def\H{{\Bbb H}}
\def\R{{\Bbb R}}
\def\Z{{\Bbb Z}}
\def\N{{\Bbb N}}
\def\frac#1#2{{#1\over#2}}

\def \prodscal#1,#2>{\langle #1,#2 \rangle}

\def\n#1{||#1||}

\def\bigv#1{\bigl|#1 \bigr|}

%
%



\arXiv 

\anglais
\tenpoint

\TITR THE MACROSCOPIC SPECTRUM OF NILMANIFOLDS WITH AN EMPHASIS ON THE HEISENBERG GROUPS|


\auteur Constantin Vernicos*|

\vfootnote*{\sevenrm Partially supported by european project ACR OFES number 00.0349 and grant of the FNRS 20-65060.01}



\centerline{\bf Abstract}
\medskip
\centerline{\vbox{\hsize 10cm \baselineskip=2.5ex \hautdepage Take a riemanniann nilmanifold, lift
its metric on its universal cover. In that way one obtains a metric invariant 
under the action of some co-compact subgroup. We use
it to define metric balls and then study the spectrum of the laplacian
for the dirichlet problem on them. 
We describe the asymptotic behaviour of the
spectrum when the radius of these balls goes to infinity. Furthermore
we show that the first macroscopic eigenvalue is bounded from above,
by an uniform constant for the three dimensional heisenberg group, and
by a constant depending on the Albanese's torus for the other nilmanifolds.
We also show that the Heisenberg groups belong to a family of nilmanifolds, where 
the equality characterizes some pseudo left invariant metrics.}}

\vfootnote{}{{\sevensl Key-words : }{\hautdepage Spectrum of the laplacian, Nilmanifolds, Homogenisation,
Stable norm, Asymptotic volume, Albanese metric, rigidity.}}
  
\vfootnote{}{{\sevensl Classification : }{\hautdepage 53C24, 58C40, 74Q99}}



\partie{Introduction and claims}\label{part1}

\paragraphe In this article we are investigating Riemannian nilmanifolds, these are compact manifolds
obtained by taking the quotient of a nilpotent Lie group by one of its subgroups,
endowed with a riemannian metric. 

Our aim is to find as much informations as we can by just looking at the balls of great radius
on the universal covering. One could believe that it is an object too shrewd to explore,
but let us recall the theorem due to Brooks \cite{brooks} (see also Sunada \cite{sunada}),
which states that if the bottom of the spectrum
of the Laplacian acting on functions on the universal cover
of a compact manifold is zero then the fundamental group is amenable. 

One could transform the statement by just saying that if the bottom of the
spectrum on the balls goes to zero as their radius goes to infinity, then
the fundamental group is amenable. 

\proclaim Question 1. Can one precise the speed of convergence to the bottom
of the spectrum on the universal cover with
respect to the radius and can one extract more geometric informations ?

The first step with that question in mind is to try to extract more
information feeding the problem with some more geometric assumptions. As
the nilpotent groups are among the amenable groups the simplest one
it is thus logical to try and explore that case first. 
In this article we answer question $1$ in the following way :

\proclaim Theorem \theo.
Let $(M^n,g)$ be nilmanifold, $B_g(\rho)$ the induced Riemannian
ball of radius $\rho$
on its universal cover and $\lambda_1\bigl(B_g(\rho)\bigr)$ the first
eigenvalue of the Laplacian on $B_g(\rho)$ for the Dirichlet problem.
\endgraf
Then
\item{1. } $\lim_{\rho\to +\infty} \rho^2\lambda_1\bigl(B_g(\rho)\bigr) = \lambda_1^\infty \leq \lambda_1(g,Alb)$
\item{2. } in case of equality the stable norm coincides with the
albanese metric.
\endgraf
Where $\lambda_1(g,Alb)$ is the first eigenvalue of the Kohn laplacian arising from
the Albanese metric on the unit ball of the carnot-caratheodory ball arising
from the same metric. Furthermore for tori and the $3$-dimensional
Heisenberg group this is a constant independent of the metric $g$.\label{theoup}

Thus we are naturally lead to explore the equality case (see \cite{ver1} for the case
of tori). Feeding with new assumptions on the
nilmanifold we get for example the following theorem (see also section \ref{example}),
introducing a family of metric which are not far from being left invariant
(see definition \ref{defpli}). 

\proclaim Theorem \theo. 
For any $2$-step nilmanifold whose center is one dimensional, the albanese metric and
the stable norm coincides if and only if the metric is pseudo left invariant.\label{thealst}

In fact it appears that we can say as much for all the eigenvalues, and the asymptotic 
behaviour of the eigenvalue in
theorem \ref{theoup} is just a particular case of the following theorem :

\proclaim Theorem \theo.
Let $(M^n,g)$ be nilmanifold, $B_g(\rho)$ the induced Riemannian
ball of radius $\rho$
on its universal cover and $\lambda_i\bigl(B_g(\rho)\bigr)$ the $i^{\rm th}$
eigenvalue of the Laplacian on $B_g(\rho)$ for the Dirichlet problem.
\endgraf
Then there exists an hypoelliptic operator $\Delta_\infty$ (the Kohn Laplacian
of a left invariant metric), whose
$i^{\rm th}$ eigenvalue for the Dirichlet problem on stable ball
is $\lambda_i^\infty$ and such that \label{theo1}
$$
\lim_{\rho\to \infty} \rho^2 \lambda_i\bigl(B_g(\rho)\bigr)= \lambda_i^\infty
$$

\paragraphe Instead of studying the Laplacians we could just study the
volume of the balls on the universal covering. As surprising at it mays
seems this also gives important informations, for example a
theorem of M.~Gromov \cite{gromov} states that
if the growth of the geodesic balls on the universal covering of a compact
manifold is polynomial, then the manifold is almost nilpotent.

\proclaim Question 2. Can we describe more precisely the asymptotic behaviour
of the volume of balls with respect to their radius and can we extract more
geometric informations ?

For the asymptotic behaviour of the polynomial case, i.e. the
nilpotent case, the answer is given by P.~Pansu \cite{pansu} who described 
precisely the growth of the geodesic balls on the universal cover of nilmanifolds in
the following way

\proclaim Theorem [Pansu]. Let $(M^n,g)$ be a nilmanifold, $B_g(\rho)$ 
the induced Riemannian ball of radius $\rho$
on its universal cover then if $\mu_g$ is the volume and $d_h$ the
homogeneous dimension of $M^n$,  then the asymptotic volume of $g$ is
$$
{\rm Asvol}(g)=\lim_{\rho\to \infty} \frac{\mu_g\bigl(B_g(\rho)\bigr)}{\rho^{d_h}} = \mu_g(M^n) \frac{\mu\bigl(B_\infty(1)\bigr)}{\mu(D_M)}
$$
where $\mu$ is a Haar measure, $B_\infty(1)$ the stable norm's unit ball and $D_M$
a fundamental domain on the universal cover.

In the particular case of tori D.~Burago and S.~Ivanov \cite{bi2} showed that the asymptotic
volume is bounded from below by the constant arising from the flat cases, which
depends only on the dimension, and furthermore that the equality caracterises flat tori. 
In \cite{ver1} there is an alternate proof
using the macroscopic spectra 
in the $2$-dimensional case.

In this article we also give a lower bound on the asymptotic volume, but
even if the lower bound is not uniform (and can not anyway, because except for
tori we can find metrics such that the asymptotic volume is as small as we want) we still
have some information on the equality case.

\proclaim Theorem \theo. 
Let $(M^n,g)$ be a nilmanifold, then its asymptotic
volume satisfies the following :
\item{1. } ${\rm Asvol (g)} \geq \displaystyle{{\rm Vol_g}(M^n){ \mu_2\bigl(B_2(1)\bigr) \over \mu_2(D_M)}}$
\item{2. } In case of equality the stable norm coincides with the albanese metric. 
\endgraf
where $\mu_2$ is the (sub-riemannian) 
measure associated to the sub-riemannian distance obtained
from the albanese metric of $M^n$ on the universal cover of $M^n$,
$B_2(1)$ is the unit ball of that distance and $D_M$ is a fundamental domain
of the universal cover of $M^n$.\label{theolast}

The last remark is that in the light of the case of tori, this should not
be the best inequality we can expect.
However even if we find the best one,
 the case of equality might not characterize left invariant metrics,
but only a particular family of metrics (see \cite{pac} for the intuition), just
like in the case of the macroscopic spectrum.


\partie{General facts about nilmanifolds\\ and their geometry}\label{part2}

Our main object of study are Nilmanifolds. In all this article a nilmanifold
will be a riemanniann compact manifolds $(M^n,g)$ which is obtained by taking the
quotient of a Lie group $G$, whose Lie Algebra is nilpotent, by one
of its co-compact subgroups $\Gamma$. We would like to stress that unlike in
other work, we don't put any restriction on the metric $g$, i.e. the metric
need not be invariant by the action of $G$ on the left (the whole point being
to characterize those metrics among all the metric one can put on $M$).

Notice that following this definition the universal covering of $M^n$ is $G$.
We will let $\tilde g$ be the metric $g$ lifted on $G$.

\subpartie{Nilpotent Lie Algebras}
\paragraphe Let  $\mathfrak{u}$ be a Lie algebra, one says that it is nilpotent if 
the sequence defined by
$$
 \mathfrak{u}^1=\mathfrak{u}, \qquad
\mathfrak{u}^{i+1}=[\mathfrak{u}^i,\mathfrak{u}].
$$ 
is such that for some $r\in \N$, $\mathfrak{u}^{r+1}=\{0\}$.

\paragraphe Among the nilpotent Lie algebra, there is a family which distinguish itself,
the graded nilpotent Lie algebras, these are
the algebras $\mathfrak{u} $ with the following decomposition :\label{graded}
$$
\mathfrak{u}=V_1\oplus\dots\oplus V_r, 
$$
such that 
\item{1--} $V_i$ is a complement of $\mathfrak{u}^{i+1}$ in $\mathfrak{u}^i$ ;
\item{2--} $[V_i,V_j]\subset V_{i+j}$ ;

\paragraphe What is quite important in our work is the fact that 
to such a graduation one can attach a one group of automorphisms $(\tau_\rho)_\rho$
called "Dilatations" such that :
$$
{\tilde \tau}_\rho(x)=\rho^ix \quad {\rm for \ all\ } x\in V_i.
$$
In fact the existence of such a family of dilatations 
is equivalent to the existence of a graduation. These dilatations plays
the same role than the dilatation in the Euclidean space.

\paragraphe All nilpotent Lie algebras are not graded. But to each nilpotent Lie
algebras we can associate a graded nilpotent one in the following way :\label{dimhomo}
$$
\mathfrak{u}_\infty = \sum_{i\geq1} \mathfrak{u}_i \slash \mathfrak{u}_{i+1}
$$
the Lie bracket being induced.
We will write
$$
\tilde \pi : \mathfrak{u} \to \mathfrak{u}_\infty
$$
and we will call {\sl Homogeneous dimension of} $\mathfrak{u}$ the number
$$
d_h= \sum_i i \dim \bigl( \mathfrak{u}^i \slash \mathfrak{u}^{i+1} \bigr) 
$$ 
\paragraphe There is another way to make that graded Lie algebra appear. \label{baseX}
Let us take a nilpotent Lie algebra $\mathfrak{u}$, remark
that for all $i$, $\mathfrak{u}^{i+1} \subset \mathfrak{u}^i$, and take independent vectors 
$X_{d_1+\dots+d_{i-1}+1},\dots,X_{d_1+\dots+d_{i-1}+d_i}$
such that the vector space $V_i$  
they span is the complement of $\mathfrak{u}^{i+1}$ in $\mathfrak{u}^i$. 
In that way one gets a basis $(X_i)$ of $\mathfrak{u}$. 
Now we define an application $\tau_\rho$ by
$$
\tau_\rho(X_p) = \rho^{\alpha(p)} X_p
$$
where $\alpha(p) = i$ if $ d_{i-1}< p \leq d_i$ with  $d_0 =0$. 

\paragraphe We also define a new Lie algebra $\mathfrak{u}_\rho$ by changing the Lie bracket in the
following way : for any $X$ and $Y$ in $\mathfrak{u}_\rho$, $[X,Y]_\rho = \tau_{1/\rho}[\tau_\rho X,\tau_\rho Y]$.
Thus $\tau_\rho$ becomes a Lie algebra isomorphism from 
$\mathfrak{u}_\rho=\bigl( \mathfrak{u}, [\cdot,\cdot]_\rho\bigr)$  to  $\bigl( \mathfrak{u}, [\cdot,\cdot] \bigr)$.

Now as $\rho$ goes to Infinity $\mathfrak{u}_\rho$ goes to $\mathfrak{u}_\infty$ in the sense that
for $i,j = 1,\dots,n$ we have
$$
[X_i,X_j]_\infty={\rm pr}_{V_{\alpha(i)+\alpha(j)}} [X_i,X_j]
$$

Notice that the graded Lie algebra is the same for any member of this family. We will
write $\tilde \pi_\rho$ the projection from $\mathfrak{u}_\rho$ to $\mathfrak{u}_\infty$ (in fact
we could avoid the subscript in $\tilde \pi_\rho$, because we can identify the Lie algebras as
linear spaces)

\paragraphe 
Notice that if the Lie algebra is graded then
$[X,Y]_\rho=[X,Y]$ and $\tau_\rho$ is a Lie algebra automorphism. Otherwise remark
that for all $X\in \mathfrak{u}_\rho$ \label{deltaretpi}
$$
\tilde \pi\bigl(\tau_\rho(X)\bigr) = {\tilde \tau_\rho}\bigl(\tilde \pi_\rho(X)\bigr)
$$

\subpartie{Remarks on exponential coordinates}
\paragraphe Let $G$ be the Lie group associated with the nilpotent Lie algebra $\mathfrak{u}$.
Thanks to the exponential coordinates we can identify $G$ with some $\R^n$,  
as a differential manifold :
$$
\phi:\R^n \to G,\ \phi:x=(x_1,\dots,x_n) \mapsto \exp x_1X_1\dots \exp x_nX_n
$$
because for nilpotent Lie groups, the exponential is a diffeomorphism between the Lie
algebra and the Lie group. Let $\ln$ be the inverse and $X_i^*$ the dual form of $X_i$

\paragraphe Moreover, if we denote by $\delta_\rho$ the following family of dilatations \label{terreaux}
$$
\delta_\rho(x_1,\dots,x_n)=(\rho^{\alpha(1)}x_1,\dots,\rho^{\alpha(n)})
$$
and we define a family of group products $*_\rho$ by setting
$$
x*_\rho y=\delta_{1/\rho}[\delta_\rho(x)\delta_\rho(y)]
$$
and also
$$
x*_\infty y=\lim_{\rho\to \infty} x*_\rho y
$$
we get a family of nilpotent Lie groups $G_\rho=(G,*_\rho)$, $\rho\in \overline\R$, whose associated
family of Lie algebras are isomorphic to the family $\mathfrak{u}_\rho$. We also
write $\pi_\rho: G_\rho\to G_\infty$ the application which sends $x\in G_\rho$ to $x \in G_\infty$. Notice
also that $d\delta_\rho=\tau_\rho$.

\paragraphe Let us define $\varphi_i:G\mapsto\R$ by $\varphi_i(g)= X_i^*\ln(g)$, then using the
Campbell-Haussdorff formula,
i.e. 
$$
\ln(x*y)=\ln(x)+\ln(y) + \frac{1}{2}\bigl[\ln(x),\ln(y)\bigr]+ C(x,y)
$$
where $C(x,y) \in \mathfrak{u}^3$, for $i=1,\dots,d_1$ we get : 
$$
X_i^*\ln(x*y) =X_i^*\ln(x)+ X_i^*\ln(y)
$$
which means that for $i=1,\dots,d_1$ the function $\varphi_i$ is a group morphism  hence
$d\varphi_i$ is left invariant, i.e. ${d\varphi_i}_{|\gamma*g}\cdot{dl_\gamma}_{|g} ={d\varphi_j}_{|g}$ thus   
$$
X_i\!\cdot\!\varphi_j= {d\varphi_j}_{|g} \!\cdot\! X_i(g)= {d\varphi_j}_{|e}\! \cdot\! X_i(e) = X_j^*\!\cdot\!X_i =\delta_{ij}
$$

\paragraphe In other words, if we use exponential coordinates, 
taking as a maximum familly of independent vector in 
$V_1=\mathfrak{u}^1 \backslash \mathfrak{u}^2$ the vectors
$X_1,\dots,X_{d_1}$, the previous calculation says that for $j=1,\dots, d_1$
and all $i$, 
$$
X_i\cdot x_j = \delta_{ij}
$$

\subpartie{Horizontal distribution and the Stable norm}

\paragraphe On the graded nilpotent Lie group $G_\infty$ 
associated to $G$ we obtain a natural distribution 
by left multiplication of
$V_1=\mathfrak{u}_1 \slash \mathfrak{u}_2\subset \mathfrak{u}_\infty$, 
we shall call that distribution {\sl Horizontal} and
write it ${\cal H}$. \label{horizontald}

\paragraphe Let us remark that because of the nilpotency and the graduation 
of the Lie algebra $\mathfrak{u}_\infty$,
a base of $V_1$ satisfies the so called H\"ormander conditions in the
Lie group $G_\infty$. 

\paragraphe We will say that a function $f$ is periodic with respect
to $\Gamma$ (the cocompact subgroup) if for every $\gamma\in \Gamma$ and $x\in G$ we have $f(\gamma*x)=f(x)$. 
Thus $\tilde g$ is periodic with respect to $\Gamma$.\label{periodicwrtg}

\paragraphe We recall what the {\sl stable norm} is :

\proclaim Definition \theo. 
Let $\n{\cdot}_\infty^*$ be the quotient of the sup norm on $1$-forms, arising from
the metric $g$, on the cohomology $H^1(M^n,\R)$. Then its dual norm on the homology $H_1(M^n,\R)$,
is called the stable norm.

\paragraphe By a theorem of K.~Nomizu \cite{nomizu}, $H_1(M^n,\R) \equiv V_1$, hence
we will call stable ball $B_\infty(1)$, the unit ball for the carnot-caratheodory
metric $d_\infty$ induced by the stable norm of $(M^n,g)$ on $G_\infty$. \label{binfini}

\subpartie{The horizontal distribution and the metric}

\paragraphe We will need to consider the family $\{G_\rho\}_{\rho\in \R}$ of simple connected Lie groups 
who corresponds to the family $\{\mathfrak{u}_\rho\}_\rho$ of Lie algebra associated to $\mathfrak{u}$.
On each $G_\rho$ we are going to pull back the metric of $G$, $\tilde g$  
and rescale it in the following way
$$
g_\rho = {1\over \rho^2} (\delta_\rho)^*{\tilde g}
$$
thus we are also able to focus on the riemanniann spaces $(G_\rho,g_\rho)$. The Laplacian
in that space will be $\Delta_\rho$.

\paragraphe If $e\in G$ is the unit element and $X \in \mathfrak{u}$ then for $\rho\in \overline \R$, $X^\rho$ will
be the $*_\rho$ left invariant field in $G_\rho$ such that $X^\rho(e)=X(e)$. Thus to the
base $(X_i)$ defined in \ref{baseX} we will associate the $*_\rho$ left invariant fields
$(X_i^\rho)$. Notice also that
$$
d\delta_\rho(X^\rho_i)=\tau_\rho(X^\rho_i)=\rho^{\alpha(i)} X_i
$$

\paragraphe Let us write the metric $\tilde g$ in the basis $(X_i)$.
$$
\tilde g=(g_{ij})
$$
We can distinguish two distincts
parts when writing the laplacian in coordinates 
$$
- \det \tilde g \Delta = \sum_{1\leq i,j \leq d_1} X_i \bigl(\det g\cdot g^{ij} X_j\bigr) + %
\sum_{\alpha(i)+\alpha(j)>2}\bigl(X_i \det g\cdot g^{ij} X_j\bigr) 
$$
which will be of significant importance in what follows.

\paragraphe It is a straightforward calculation to find that the metric $g_\rho$ in the 
coordinates $(X_i^\rho)$ is written (for $x\in G_\rho$ and $\rho\in \R$):
$$\eqalign{%
-\det \tilde g(\delta_\rho x) \Delta_\rho = &\sum_{1\leq i,j \leq d_1} X_i^\rho \bigl(\det g~g^{ij} (\delta_\rho x)  X_j^\rho\bigr)+\cr%
&\sum_{\alpha(i)+\alpha(j)>2}\rho^{2-\alpha(i)-\alpha(j)}X_i^\rho\bigl( \det \tilde g~g^{ij}(\delta_\rho x)X_j^\rho\bigr)\cr} \numeq \label{fsublap} 
$$
In this formula the whole difference between the two parts becomes clear, indeed it is
quite apparent now that the second part vanishes when $\rho$ goes to infinity.

\paragraphe That is the reason why we introduce $\nabla_{\cal H}$ by :
$$
\nabla_{\cal H} f = (X_1^\infty\cdot f,\dots,X_{d_1}^\infty\cdot f)
$$

\subpartie{Gromov-Haussdorff convergence of balls}\label{ghcball}

\paragraphe In what follows $B_g(\rho)$ will always mean the geodesic ball of radius $\rho$ on
the universal cover of $(M^n=\Gamma\backslash G,g)$, and $B_\rho(1)$ will be the geodesic ball of radius one
on $(G_\rho,g_\rho)$. We will also need to define $\mu_\rho$ (resp. $\mu_g$) the Riemannian volume associated to
$g_\rho$ (resp. $g$) and $\mu_\infty$ defined as follows : Let $D_\Gamma$ be a fundamental domain in $G$
and $\mu$ a Haar measure on $G_\infty$ then
$$
\mu_\infty=\frac{\mu_g(D_\Gamma)}{\mu\bigl(\pi(D_\Gamma)\bigr)}\ \mu
$$

\proclaim Theorem \theo. Let $(M^n,g)$ be a riemanniann nilmanifold,
$d_g$ the induced distance on its universal cover and $(\delta_\rho)_{\rho\in\R}$ the familly of
isomorphisme between $G_\rho$ and $G$ if one writes for any $x,y \in G_\rho$ : 
$$
d_\rho(x,y)= {d_g(\delta_\rho x,\delta_\rho y)\over \rho} 
$$
and $B_\rho(1)$ the unit ball for each of this rescaled distances,
then the family of metric spaces $(B_\rho(1),d_\rho,\mu_\rho)$ converges in the Gromov-Haussdorff measure
topology to $(B_\infty(1),d_\infty,\mu_\infty)$ as $\rho$ goes to infinity.

\proof
The Gromov-Haussdorff convergence comes from P.~Pansu work \cite{pansu}, which implies
$$
\lim_{\rho\to \infty} \frac{d_\infty\bigl(\pi_\rho(x),\pi_\rho(y)\bigr)}{d_\rho(x,y)} = %
\lim_{\rho\to \infty} \frac{d_\infty\bigl(\pi\circ\delta_\rho (x), \pi\circ\delta_\rho(y) \bigr)}{d_g(\delta_\rho x,\delta_\rho y)}=1
$$
It remain to show the measure part.

\proclaim Claim.
for any compact domain $A$ in $G_\infty$, whose boundary is of haar measure $0$, and any
function $f \in L^1(A,\mu_\infty)$ we have
$$
\lim_{\rho\to \infty} \int_{\pi_\rho^{-1}(A)} f\bigl( \pi_\rho(x)\bigr) d\mu_\rho(x) = \int_A f d\mu_\infty 
$$

Indeed, let $A$ be a domain in $G_\infty$, then $\pi_\rho^{-1}(A)$ belongs to $G_\rho$ and $\delta_\rho\circ \pi_\rho^{-1}(A)$ belongs
to $G$. We will write $*$ the law group of $G$. Let $z_1,\dots,z_k$ and $\zeta_1,\dots,\zeta_l$
be elements of $\Gamma$ such that $\zeta_j*D_\Gamma \cap \delta_\rho\circ \pi_\rho^{-1}(A) \not= \emptyset$ for any $j$ and
$$
\bigcup_i z_i*D_\Gamma \subset \delta_\rho\circ \pi_\rho^{-1}(A) \subset  \bigcup_j \zeta_j*D_\Gamma
$$
Let us notice that
$$
\mu_g(D_\Gamma)=\frac{\mu_g(D_\Gamma)}{\mu\bigl(\pi(D_\Gamma)\bigr)}\mu\bigl(\pi(D_\Gamma)\bigr) = \mu_\infty\bigl(\pi(D_\Gamma)\bigr)
$$
then we get
$$\displaylines{%
\sum_i \inf_{\delta_\rho\circ \pi_\rho^{-1}(x)\in z_i*D_\Gamma} f(x) \mu_\infty\bigl(\pi(D_\Gamma)\bigr) \leq \hfill \cr
\hfill \int_{\delta_\rho\circ \pi_\rho^{-1}(A)} f\bigl({\tilde \delta}_{1/\rho}\circ\pi(x)\bigr) d\mu_g(x) \leq %
\sum_j \sup_{\delta_\rho\circ \pi_\rho^{-1}(x)\in \zeta_j*D_\Gamma} f(x) \mu_\infty\bigl(\pi(D_\Gamma)\bigr)}
$$
divide all members by $\rho^{d_h}$ (see \ref{dimhomo}) we get : 
$$\displaylines{%
\sum_i \inf_{x  \in \pi_\rho\circ\delta_{1/\rho}(z_i*D_\Gamma)} f(x) \mu_\infty\bigl({\tilde \delta}_{1/\rho}\circ\pi(D_\Gamma)\bigr) \leq \hfill \cr
\hfill \int_{\pi_\rho^{-1}(A)} f\bigl(\pi_\rho(x)\bigr)d\mu_\rho \leq %
\sum_j \sup_{x \in  \pi_\rho\circ\delta_{1/\rho}(\zeta_j*D_\Gamma)}  f(x) \mu_\infty\bigl({\tilde \delta}_{1/\rho}\circ\pi(D_\Gamma)\bigr)}
$$
then the extremal terms are riemann's sums which converge toward
$\int_A f d\mu_\infty \ {\rm .}$

Using the claim and the fact that the functionals $x\mapsto d_\rho(0,\pi_\rho^{-1}(x))$ converges
simply toward the functional $x\mapsto d_\infty(0,x)$ on $B_\infty(1)\backslash \partial B_\infty(1)$ we can conclude. 
\qed

\paragraphe Let us introduce the asymptotic volume as
$$
{\rm Asvol}(g) = \lim_{\rho\to \infty} \frac{\mu_g\bigl(B_g(\rho)\bigr)}{\rho^{d_h}}
$$
following our last statement it is straightforward that
$$
{\rm Asvol}(g) = \mu_\infty\bigl(B_\infty(1)\bigr)
$$
\partie{Looking for convergences}

For the proof of theorem \ref{theo1} we will a few more definitions. Indeed
we are going to look at a family of $L^2$-spaces which are not defined on
the same space. More exactly we are going to look at $L^2(B_\rho(1),\mu_\rho)$. However
we are going to show that in some sens (see below \ref{ghcball}) these
spaces converge towards $L^2(B_\infty(1),\mu_\infty)$. So now we would like to give
a precise meaning to the fact that a net composed of functions in $L^2(B_\rho(1),\mu_\rho)$
converges toward a function in $L^2(B_\infty(1),\mu_\infty)$.

Once this will be done, we will be able to give a meaning to the convergence
of a net of operators (the resolvent of the Laplacian for example).

\subpartie{Convergence on a net of Hilbert spaces}

In what follows ${\cal A}$ and ${\cal B}$ are directed set.

\paragraphe Let $(X_\alpha,d_\alpha,m_\alpha)_{\alpha\in {\cal A}}$ be a net of compact measured metric spaces converging in
the Gromov-Hausdorff measured topology to $(X_\infty,d_\infty,m_\infty)$. We will write $L_\alpha^2=L^2(X_\alpha,m_\alpha)$
(resp. $L^2_\infty(X_\infty,m_\infty)$) for the square integrable function spaces. Their respective scalar product
will be $\langle\cdot,\cdot\rangle_\alpha$ (resp.  $\langle\cdot,\cdot\rangle_\infty$) and $\n{\cdot}_\alpha$ (resp. $\n{\cdot}_\infty$). 

\paragraphe Furthermore we suppose that in every  $L_\alpha^2$ the continuous functions  
form a dense subset $C^0(X_\alpha)$.

\proclaim Definition \theo. We say that a net $(u_\alpha)_{\alpha \in {\cal A}}$ of functions $u_\alpha\in L^2_\alpha$ strongly
converges to $u\in L^2_\infty$ if there exists a net $(v_\beta)_{\beta\in {\cal B}} \subset C^0(X_\infty)$ converging to $u$ in $L^2_\infty$
such that \label{deficonvfort}
$$
\lim_\beta \limsup_\alpha \n{f^*_\alpha v_\beta -u_\alpha}_\alpha =0;
$$
where $(f_\alpha)$ is the net of Hausdorff approximations. We will also talk of strong convergence in 
${\cal L}^2$.

\proclaim Definition \theo. We say that a net $(u_\alpha)_{\alpha \in {\cal A}}$ of functions $u_\alpha\in L^2_\alpha$ weakly
converges to $u\in L^2_\infty$ if and only if for every net $(v_\alpha)_{\alpha\in {\cal A}}$ strongly converging to $v\in L^2_\infty$
we have \label{deficonvfaib}
$$
\lim_\alpha \langle u_\alpha,v_\alpha\rangle_\alpha =\langle u,v\rangle_\infty \numeq
$$
We will also talk of weak convergence in ${\cal L}^2$.

The following claim (whose proof we don't give, see \cite{vernicos} pages 32---33) justifies those two definitions

\proclaim Lemma \theo.
Let $(u_\alpha)_{\alpha \in {\cal A}}$ be a net of functions $u_\alpha\in L^2_\alpha$. If $(\n{u_\alpha}_\alpha)$ is
uniformly bounded, then there exists a weakly converging subnet.
Furthermore every weakly converging net is uniformly bounded.\label{subwcvnet}

\paragraphe Now that we gave sense to the convergence of a net of functions, we
are going to define convergences of a net of operators. 
Let $B_\infty\in {\cal L}(L^2_\infty)$ and $B_\alpha\in {\cal L}(L^2_\alpha)$ for every $\alpha\in {\cal A}$.
\proclaim Theorem and Definition \theo. Let $u,v \in L^2_\infty$ and $(u_\alpha)_{\alpha\in{\cal A}}$, $(v_\alpha)_{\alpha\in{\cal A}}$ two nets
such that $u_\alpha, v_\alpha\in L^2_\alpha$. We say that the net of operators $(B_\alpha)_{\alpha\in{\cal A}}$ strongly 
(resp. weakly, compactly) converges to $B$ if $B_\alpha u_\alpha \to Bu$ strongly (resp. weakly, strongly) for
every net $(u_\alpha)$ strongly (resp. weakly, weakly) converging to $u \iff$
$$
\lim_\alpha \langle B_\alpha u_\alpha,v_\alpha\rangle_\alpha = \langle Bu,v\rangle_\infty \numeq
$$
for every $(u_\alpha)$, $(v_\alpha)$, $u$ and $v$ such that $u_\alpha\to u$ strongly (resp. weakly, weakly) and
$v_\alpha\to v$ weakly (resp. strongly, weakly), (See \cite{vernicos} page 35 for the justifications).

\subpartie{The importance of being compactly convergent}

The theorem and definition of the previous section is useful, 
for the goal we would like to achieve thanks
to the following one (whose proof can be found in  \cite{ver1} for example), which links
the convergence of the resolvents $R^\alpha_\zeta$ associated to a family of functionals $A_\alpha$, with their spectra.

\proclaim Theorem \theo. Let $R^\alpha_\zeta\to R_\zeta$ compactly for all $\zeta$ outside
the spectra of $(A_\alpha)$. Assume that all
resolvents $R^\alpha_\zeta$ are compact. Let $\lambda_k$ (resp. $\lambda_k^\alpha$) be the $k^{\rm th}$
eigenvalue of $A$ (resp. $A_\alpha$) with multiplicity. We take $\lambda_k=+\infty$ if
$k>\dim L^2_\infty +1$ when $\dim L^2_\infty< \infty$ and $\lambda_k^\alpha=+\infty$ if
$k>\dim L^2_\alpha +1$ when $\dim L^2_\alpha< \infty$. Then for every $k$
$$
\lim_\alpha \lambda_k^\alpha=\lambda_k
$$
Furthermore let $\{\varphi_k^\alpha\mid k=1,\ldots,\dim L_\alpha^2\}$ be an orthonormal bases of $L^2_\alpha$ such that
$\varphi_k^\alpha$ is an eigenfunction of $A_\alpha$ for $\lambda_k^\alpha$. Then there is a sub-net such that
for all $k\leq \dim L^2_\infty$ the net $(\varphi_k^\alpha)_\alpha$ strongly converges to the eigenfunction 
$\varphi_k$ of $A$ for the eigenvalue $\lambda_k$, and such that the family
 $\{\varphi_k\mid k=1,\ldots,\dim L_\alpha^2\}$ is an orthonormal basis of $L^2_\infty$.\label{cvssets}

\subpartie{What shall we finally study ?} \label{respell}

We are now going to focus on the spectrum of the balls $B_g(\rho)$, and we
want to show that the eigenvalues of the laplacian are converging to
zero with a $1/\rho^2$ speed, to be more specific we want to find
a precise equivalent.

For this let recall that $\Delta_\rho$ is the Laplacian (or Laplace-Beltrami operator) 
associated to the rescaled metrics $g_\rho=1/\rho^2(\delta_\rho)^*g$ on $G_\rho$,
and for any function $f$ from $B_g(\rho)$ to $\R$ lets associate a function $f_\rho$ on
$B_\rho(1)$ by $f_\rho(x)=f(\delta_\rho\!\cdot\! x)$. Then it is an easy calculation to see that for any $x\in B_\rho(1)$ :
$$
\rho^2\bigl(\Delta f\bigr)(\delta_\rho\!\cdot\! x)= \bigl(\Delta_\rho f_\rho\bigr)(x)
$$
hence the eigenvalues of $\Delta_\rho$ on $B_\rho(1)$ are exactly the eigenvalues of $\Delta$ on $B_g(\rho)$
multiplied by $\rho^2$ and
our problems becomes the study of the spectrum of the laplacian
$\Delta_\rho$ on $B_\rho(1)$. 

Enlightened by what happens on tori we would like to show that
there is some operator $\Delta_\infty$ acting on $B_\infty(1)$ (see \ref{binfini}) such that, in a good sense, the net
of laplacian $(\Delta_\rho)$ converges towards $\Delta_\infty$ and the spectra also
converge to the spectrum of $\Delta_\infty$. 

In the light of theorem \ref{cvssets}, the good sense is the compact convergence
of the resolvent like in our paper on the macroscopical
sound of tori \cite{ver1}. The proof is quite similar, but needs some
adaptation to the geometry of nilmanifolds. The following
section gives the proof.

\partie{Homogenisation and proof of theorem \ref{theo1}}\label{proof1}

The first step consisted in showing the convergence of
the metric geodesic balls w.r.t. the Gromov-Haussdorff measure topology (see \ref{ghcball}).

The next step will consist in 
introducing some functions linked with the Albanese metric (see \ref{homogelap}) which
will lead us to the definitions of the functional spaces involved and the compact
inclusion we can deduce (see \ref{asympcomp}). Finaly we prove the compact
convergence of the resolvent (see \ref{compconv}), which thanks to theorem~\ref{cvssets} finishes
the proof.

\subpartie{Homogenisation of the Laplacian and Albanese's Torus }\label{homogelap} 
In this section we are going to "built" the operator $\Delta_\infty$ of theorem \ref{theo1}
and give the proof of theorem \ref{theolast}. 

\paragraphe Let $D_\Gamma$ be  a fundamental domain. 
Let $\chi^i$ be the unique periodic with respect to $\Gamma$ solution 
(up to an additive constant) of (for $1\leq i\leq r$)  
$$
\Delta\chi^i = \Delta x_i \ {\rm on}\ D_\Gamma
$$
The operator $\Delta_\infty$ is then defined by (we use Einstein's summation convention)
$$
\Delta_\infty f = -{1 \over {\rm Vol}(g)} \sum_{1\leq i,j\leq d_1} \Biggl(\int_{D_\Gamma} g^{ij}-g^{ik} X_k\cdot\chi^j~d\mu_g \Biggr)%
{X_i^\infty\cdot X_j^\infty f} \numeq
$$ 
Now let us write $\eta_j(x)=\chi^j(x)-x_j$ the induced harmonic function and
$$
q^{ij}= {1 \over {\rm Vol}(g)} \Biggl(\int_{D_\Gamma} g^{ij}-g^{ik} X_k\cdot\chi^j~d\mu_g \Biggr)
$$
we can notice that the $d\eta_i$ are harmonic $1$-forms on the nilmanifold.
It is not difficult now to show that

\proclaim Proposition \theo.
Let $\langle\cdot,\cdot\rangle_2$ be the scalar product induced on $1$-forms by the Riemannian
metric $g$. Then 
$$
q^{ij}= {1\over {\rm Vol}(g) } \langle d\eta_i,d\eta_j\rangle_2 =q^{ji} 
$$
thus $\Delta_\infty$ is an Hyppoelliptic operator.

\paragraphe In fact we can say more, $(q^{ij})$ induces a scalar product on harmonic $1$-forms
(whose norm will be written $\n{\cdot}_2$) and then to $H^1(M^n,\R)$, which can be identified
with the horizontal ${\cal H}$ following a theorem of K.~Nomizu \cite{nomizu}. Indeed,
as mentioned earlier, we can see the $(d\eta_i)$ as $1$-forms over the nilmanifold.
Being a free family they can be seen as a basis of $H^1(M^n,\R)$ (Hodge's theorem).
Thus by duality this yields also a scalar product $(q_{ij})$ over $H_1(M^n,\R)$ 
(whose induced norm will be written $\n{\cdot}^*_2$). 

\paragraphe The question naturally arising
is to know how is this norm related to the stable norm. To understand their link we have to
go back on $H^1(M^n,\R)$. Indeed the stable norm is the dual of the norm
obtained by quotient of the sup norm on $1$-forms (see Pansu \cite{pansu2} lemma 17),
which we write $\n{\cdot}_\infty^*$, and the norm $\n{\cdot}_2$ comes from the
normalised $L^2$ norm. Thus mixing the H\"older inequality and the Hodge-de~Rham
theorem we get :
\proclaim Proposition \theo. For every $1$-form $\alpha$ we have
$$
\n{\alpha}_2 \leq \n{\alpha}_\infty^*
$$
thus by duality, for every $\gamma\in H_1(M^n,\R)$ we have \label{alinst}
$$
\n{\gamma}_\infty \leq \n{\gamma}^*_2
$$
in other words the unit ball of the sub-riemannian metric arising from
 $\n{\cdot}^*_2$ is included in $B_\infty(1)$.

\paragraphe Let us recall that the manifold $H_1(M^n,\R)/ H_1(M^n,\Z)$
with the flat metric induced by $\n{\cdot}^*_2$ is usually called the Jacobi manifold
or the Albanese torus of $(M^n,g)$. This last proposition also
implies the following inequality, regarding the asymptotic volume

\proclaim Corollary \theo. 
Let $(M^n,g)$ be a nilmanifold, then its asymptotic
volume satisfies the following inequality :
$$
{\rm Asvol (g)} \geq {\rm \mu_g}(M^n){ \mu_2\bigl(B_2(1)\bigr) \over \mu_2(D_\Gamma)} 
$$
where $\mu_2$ is the (sub-riemanniann) 
measure associated to the sub-riemanniann distance obtained
from the albanese metric of $M^n$ on the universal cover of $M^n$,
$B_2(1)$ is the unit ball of that distance and $D_\Gamma$ is a fundamental domain
of the universal cover of $M^n$.

\proof
Following Nomizu \cite{nomizu} we can identify the horizontal space ${\cal H}$
with $H_1(M^n,\R)$. This allows us to get two sub-riemannian distances
$d_2$ and $d_\infty$ from $\n{\cdot}_2^*$ and $\n{\cdot}_\infty$ respectively. The previous
proposition implies that the ball of $d_2$ is inside the stable ball.
Thus for any Haar measure $\mu$ one gets the following inequality :
$$
\mu\bigl(B_2(1)\bigr) \leq \mu\bigl(B_\infty(1)\bigr)
$$
now taking for $\mu$ the haar measure $\mu_\infty$ (see section \ref{ghcball})  
giving the asymptotic volume we can conclude. \qed

Remark that theorem \ref{theolast} is now a simple corollary of that
last corollary.

\subpartie{Asymptotic compactness}\label{asympcomp}
\paragraphe Let us now define the various functional spaces involved. For $\rho\in \overline{\R}$,
$L_\rho^2=L^2\bigl(B_\rho(1),d\mu_\rho\bigr)$ will be the space of square integrable functions
over the ball $B_\rho(1)$, which is a Hilbert space with the scalar product
$$
(u,v)_\rho=\int_{B_\rho(1)} uv~d\mu_\rho
$$
whose norm will be $|\cdot|_\rho$. Hence ${\cal L}^2$ will be the net of spaces $(L_\rho^2)$
with either the strong or weak topoly induced by the definitions \ref{deficonvfort} and
\ref{deficonvfaib}.

\paragraphe Following the usual nomenclature we will be interested in the following spaces
$$\eqalign{
& H^1_\rho\bigl(B_\rho(1)\bigr) = \biggl\{ v \Bigm|v,{X_i^\rho\cdot v} \in L^2\bigl(B_\rho(1),d\mu_\rho\bigr),\ 1\leq\alpha(i)\leq r \biggr\} \cr
\Biggl( {\rm resp.} \qquad 
& H^1_\infty\bigl(B_\infty(1)\bigr) = \biggl\{ v \Bigm|v,{X_i^\infty\cdot v} \in L^2\bigl(B_\infty(1),d\mu_\infty\bigr),\ 1\leq i\leq d_1 \biggr\} \Biggr)%
\cr}
$$
which becomes Hilbert spaces with the norm $\n{\cdot}_\rho$ defined by
$$\eqalign{
&\n{v}_\rho^2= |v|_\rho^2 + \sum_{1\leq\alpha(i)\leq r}^n \bigv{X_i^\rho\cdot v}_\rho^2 \cr
\Biggl( {\rm resp.} \qquad &\n{v}_\infty^2= |v|_\infty^2 + \sum_{1\leq i\leq d_1}^n \bigv{X_i^\infty\cdot v}_\infty^2 \Biggr)\cr}
$$

\paragraphe Hence $H^1_{\rho,0}\bigl(B_\rho(1)\bigr)$ will be the closure
of the  $C^\infty\bigl(B_\rho(1)\bigr)$ functions with compact support, in $H^1_\rho\bigl(B_\rho(1)\bigr)$
for the norm $\n{\cdot}_\rho$.

\paragraphe  We can define a "spectral structure" on 
$L^2_\rho$ by expanding the Laplacian (sub-laplacian for $\Delta_\infty$)
defined on $H^1_{\rho,0}\bigl(B_\rho(1)\bigr)$ thanks to the following quadratic form
$$
\n{v}_{\rho,0}^2 = |v|_\rho^2 + (v,\Delta_\rho v)_\rho 
$$
Now let us see what can we say of a bounded net in
$H^1_{\rho,0}\bigl(B_\rho(1)\bigr)$ for this quadratic form.

\proclaim Lemma \theo.
Let $(u_\rho)_\rho$ be a net with $u_\rho \in H^1_{\rho,0}\bigl(B_\rho(1)\bigr)$ for every $\rho$, if there is
a constant $C$ such that for all $\rho>0$ we have
$$
\n{u_\rho}_{\rho,0} \leq C
$$
then there is a strongly converging sub-net in ${\cal L}^2$.\label{lemapivot}

\proof
Let $B$ a compact set such that $\bigcup_\rho \pi_\rho\bigl(B_\rho(1)\bigr)\subset B \subset G_\infty$  we are going to show
that the strong convergence in $L^2(B,\mu_\infty)$ implies the strong convergence
in ${\cal L}^2$. Then the compact embedding of $H^1_\infty\bigl(B\bigr)$ in $L^2\bigl(B,\mu_\infty\bigr)$
will conclude the proof.

Let us first notice that the periodicity with respect to $\Gamma$, and the co-compactness
of $\Gamma$ gives the existence of two constants
$\alpha$ and $\beta$ such that (we suppose the norms defined on $B$, and identify $B$ et $\pi_\rho^{-1}B$)
$$
\alpha |v|_\infty \leq |v|_\rho \leq \beta |v|_\infty {\rm\ .}
$$
Let us start by taking a net $(v_\rho)$ strongly converging in $L^2(B,\mu_\infty)$ to $v_\infty$
we also assume $v_\rho\circ\pi_\rho \in  H^1_{\rho,0}\bigl(B_\rho(1)\bigr)$ for every $\rho$ (because it is all we need). 
The first remark to be done is that thanks to the Gromov Haussdorff convergence,
$v_\infty\in L^2_\infty$ (we mean that $v_\infty$ can be considered equal to zero outside $B_\infty(1)$).
Thus, let us take $c_p \in C_0^\infty\bigl( B_\infty(1) \bigr)$ be a sequence of 
functions strongly converging to $v_\infty$ in $L^2_\infty$. 
We have 
$$
|c_p\circ\pi_\rho-v_\rho\circ\pi_\rho|_\rho \leq \beta |c_p-v_\infty|_\infty +\gamma |v_\infty-v_\rho|_\infty
$$
now let  $\varepsilon>0$ then for $p$ large enough $\beta |c_p-v_\infty|_\infty \leq \varepsilon$. We fix $p$
large enough and take $\rho$ large enough for the second term to converge to $0$
which gives us the strong convergence needed (see \ref{deficonvfort}).

Now to conclude observe that from the assumptions the net $(u_\rho\circ \pi_\rho^{-1})$ (if need
be we extend this function by zero outside $B_\rho(1)$) is bounded in
$H^1_\infty\bigl(B\bigr)$, hence using the compact embedding of $H^1_\infty\bigl(B\bigr)$
in $L^2(B,\mu_\infty)$ (with the right regularity assumption on the boundary of $B$) 
we can extract a strongly converging net in $L^2(B,\mu_\infty)$ an by what
we just did in ${\cal L}^2$.
\qed

\subpartie{Compact convergence of the resolvents}\label{compconv}

\paragraphe Let $\lambda>0$,  $a^\rho_\lambda(u,v)=(\Delta_\rho u,v)_\rho+\lambda(u,v)_\rho$
and $G_\lambda^\rho$ be the operator from $L^2_\rho$ to $H_{\rho,0}^1\subset L^2_\rho$
such that
$$
a^\rho_\lambda(G_\lambda^\rho f,\phi) = (f,\phi)_\rho\quad \forall \phi \in H_{\rho,0}^1 {\rm.} \numeq \label{predemo1}
$$

\paragraphe We want to show that the net of operators $(G_\lambda^\rho)$ converges compactly to
$G_\lambda$ the operator corresponding to the homogenised problem:
$$
a^\infty_\lambda(G_\lambda F,\Phi)=(F,\Phi)_\infty  \quad \forall \Phi \in H_{\infty,0}^1 \numeq \label{predemo2}
$$
with $(F,\Phi)_\infty=\int_{B_\infty(1)}F\Phi~d\mu_\infty$ and
$$
a^\infty_\lambda(u,v)=\int_{B_\infty(1)} q^{ij}~X_i^\infty u~X_j^\infty~d\mu_\infty+\lambda(u,v)_\infty
$$
in other word we want to show the following theorem

\proclaim Theorem \theo.
  For every $\lambda<0$, the net of resolvents $(R_\lambda^\rho)_\rho$ of the Laplacian $(\Delta_\rho)$ converges 
compactly to $R_\lambda^\infty$, the resolvent of $\Delta_\infty$ from the homogenised problem.
Thus the net $(\Sigma_\rho)$ compactly converges to $\Sigma_\infty$.\label{cvres}

\proof
This comes from the fact that $R_\lambda^\rho=-G_{-\lambda}^\rho$ and $R_\lambda^\infty=-G_{-\lambda}$.

\noindent{\bf First step :}

Let $f_\rho$ be a weakly convergent net to $f$ in
${\cal L}^2$, lemma \ref{subwcvnet} tells us that this net is uniformly 
bounded in ${\cal L}^2$ and in
$H_\rho^{-1}$, the dual space of $H_{\rho,0}^1$.

Let $f_\rho \in H_{\rho,0}^1$ then  thanks to \ref{predemo1}  we have :
$$
\alpha \n{G^\rho_\lambda f_\rho}^2_{\rho,0} \leq (f_\rho,G_\lambda^\rho f_\rho)_\rho \leq %
K\n{f_\rho}_{H_\rho^{-1}}\n{G_\lambda^\rho f_\rho}_{\rho,0}
$$
thus
$$
\n{G_\lambda^\rho f_\rho}_{\rho,0} \leq C \n{f_\rho}_{H_\rho^{-1}}
$$
the net $(G_\lambda^\rho f_\rho)$ being uniformly bounded for the norms $\n{\cdot}_{\rho,0}$, using lemma \ref{lemapivot}
 there is a subnet strongly converging in ${\cal L}^2$.  i.e.
$$
u_\rho=G_\lambda^\rho f_\rho \to \tilde u_\lambda {\rm\ strongly\ in\ } {\cal L}^2 \numeq \label{demo1}
$$

Furthermore  $P_\rho=(g^{ij}_\rho)\nabla G_\lambda^\rho f_\rho$ is also bounded in
${\cal L}^2$ thus there is a subnet of the net $P_\rho$ 
weakly converging in ${\cal L}^2$ to $\tilde P_\lambda\in L^2_\infty$, moreover $\tilde P_\lambda$ is horizontal
for if we write $P^{i}_\rho$ and $ P^{i}_\lambda$ the coordinates of $P_\rho$ and $\tilde P_\lambda$ then we have
$P^{i}_\rho = (g^{ij}_\rho) \nabla G_\lambda^\rho f_\rho = \rho^{2-\alpha(i)-\alpha(j)} \bigl(g^{ij}(\delta_\rho x)\bigr) \nabla G_\lambda^\rho f_\rho$, so if $\alpha(i)\geq2$ then
this coordinates strongly converges to $0$ in ${\cal L}^2$, because $\bigl(g^{ij}(\delta_\rho x)\bigr) \nabla G_\lambda^\rho f_\rho$
is also bounded. 

Now for any $\phi_\infty\in L^2_\infty $ let $\phi_\rho$ be a strongly converging net to
$\phi_\infty$ in ${\cal L}^2$ then
$$
  \eqalign{%
  \int_{B_\rho(1)} P_\rho\!\cdot\!\nabla\phi_\rho~d\mu_\rho +\lambda(G_\lambda^\rho f_\rho,\phi_\rho)_\rho&= (f_\rho,\phi_\rho)_\rho \to \cr %
  \int_{B_\infty(1)} \tilde P_\lambda\!\cdot\!\nabla_{\cal H}\phi_\infty~d\mu_\infty+\lambda(u^*_\lambda,\phi_\infty)_\infty&=(f,\phi_\infty)_\infty{\rm .}\cr}\numeq \label{demo1.1}
$$
Thus it is enough to show that $\tilde P_\lambda=\bigl(q^{ij}\bigr)\nabla_{\cal H}\tilde u_\lambda$ on
$B_\infty(1)$ because it induces $\tilde u_\lambda=G_\lambda f$.

\noindent{\bf Second step : } 

We first take $\chi^k(y)$ (see \ref{homogelap}) such that ${\cal M}(\chi^k)= 0$ and we define
$$
w_\rho(x)=x_k- \frac{1}{\rho} \chi^k(\delta_\rho x) \numeq \label{demo2}
$$
for every $k=1,\dots,d_1$.
Then
$$
w_\rho \to x_k \ {\rm  strongly\ in }\  {\cal L}^2.\numeq \label{demo2.1}
$$  
and by construction  of $\chi^k$ (see \ref{homogelap})  we have
$$
 -X_i^\rho\bigl(\det \tilde g(\delta_\rho x) g^{ij}_\rho~X_j^\rho w_\rho \bigr) =0 \ {\rm on }\ B_\rho(1) {\rm .} \numeq \label{demo3}
$$
We multiply this equation by a test function 
$\phi_\rho \in V_\rho$ and after an integration we get 
$$
\int_{B_\rho(1)}g^{ij}_\rho~X_j^\rho w_\rho~X_i^\rho\phi_\rho~d\mu_\rho=0  \numeq \label{demo4}
$$ 
Let $\varphi \in C^\infty_0(B_\infty(1))$ (notice that for $\rho$ large enough 
the support of $\varphi$ will be in $\pi_\rho\bigl(B_\rho(1)\bigr)$ ) and $\phi_\rho=\varphi\circ\pi_\rho w_\rho$ which we
put into the equation \ref{predemo1} and into the equation \ref{demo4} we put $\phi_\rho=\varphi\circ\pi_\rho u_\rho$,  
and then we subtract the results
$$\displaylines{%
\qquad \int_{B_\rho(1)} g^{ij}_\rho\Bigl( {X_j^\rho u_\rho}~{\bigl(X_i^\rho(\varphi\circ\pi_\rho)\bigr)}~w_\rho - X_j^\rho w_\rho~\bigl(X_i^\rho(\varphi\circ\pi_\rho)\bigr)~u_\rho \Bigr)~d\mu_\rho \hfill \cr
\hfill =\int_{B_\rho(1)} f_\rho w_\rho\varphi\circ\pi_\rho~d\mu_\rho -\lambda \int_{B_\rho(1)}\varphi\circ\pi_\rho u_\rho w_\rho~d\mu_\rho \qquad \nume \label{demo5} }
$$

Now let  $\rho \to \infty $ in \ref{demo5}, all
terms converge because they are product of one strongly converging net and one weakly 
converging net in ${\cal L}^2$.
More precisely,
    \item{(i)} $P_\rho$, remember that $P_\rho^i=g^{ij}_\rho X_j^\rho u_\rho$, weakly converges 
  to $\tilde P_\lambda$ in ${\cal L}^2$ following \ref{demo1.1};
    \item{(ii)} $\bigl(X_i^\rho(\varphi\circ\pi_\rho)\bigr) w_\rho$ strongly converges to
  $(X_i^\infty\varphi) x_k$ in ${\cal L}^2$ from \ref{demo2.1} and because, writing $l_x^\rho$ the function
left multiplication by $x$ in $G_\rho$:
$$
X_i^\rho(\varphi\circ\pi_\rho)_{\mid x} = d\varphi_{\pi_\rho\circ l_x^\rho(e)}\circ{d\pi_\rho}_{\mid l_x^\rho(e)}\circ dl_x^\rho\cdot X_i^\rho(e)
$$
now by definition $l_x^\rho\to l_x^\infty$ and $\pi_\rho\to id_{G_\infty}$ which explains why 
$$X_i^\rho(\varphi\circ\pi_\rho) \to X_i^\infty\varphi$$ pointwise (and weakly ${\cal L}^2$ from the claim in the proof
of section \ref{ghcball}).
    \item{(iii)} for $1\leq i,j\leq d_1$, $g^{ij}_\rho X_i^\rho w_\rho$ is periodic with respect to 
$\delta_{1/\rho}\Gamma$ and weakly converges in ${\cal L}^2$ towards the mean value
  $$
  q^{jk}=\frac{1}{\mu_g(D_\Gamma)}\int_{D_\Gamma}\biggl(g^{ij}(y)\Bigl(\delta_{ik}-X_i \chi^k(y)\Bigr)\biggr)d\mu_g
  $$
this comes from the following claim
\proclaim Lemma \theo.
Let $h$ be a function periodic with respect to  $\Gamma$ on $G$. Let $h_\rho$ be defined on
$G_\rho$ by $h_\rho(x)=h(\delta_\rho x)$. Then $(h_\rho)$ weakly converges in ${\cal L}^2$ toward
$$
h_\infty=\frac{1}{\mu_g(D_\Gamma)}\int_{D_\Gamma} h d\mu_g
$$
i.e. for any $u_\rho \to u_\infty$ strongly in ${\cal L}^2$ we have
$$
\int_{B_\rho(1)} u_\rho h_\rho d\mu_\rho \to h_\infty \int_{B_\infty(1)} u_\infty d\mu_\infty
$$

To see this apply the proof in section \ref{ghcball} to $h^{1/n}g$ instead of $g$ (even
if it is not a measure, what makes everything work in that proof is the fact
that $\det g$ in the coordinates $(X_i)$ is periodic with respect to $\Gamma$).

    \item{(iv)} for $\alpha(i)+\alpha(j)>2$, $g^{ij}_\rho X_i^\rho w_\rho = \rho^{2-\alpha(i)-\alpha(j)} g^{ij}(\delta_\rho x) X_i^\rho w_\rho$ thus this
term weakly converges in ${\cal L}^2$ towards $0$.
    \item{(v)} $\bigl(X_j^\rho(\varphi\circ\pi_\rho)\bigr) u_\rho$ strongly  converges to
  $( X_j^\infty\varphi) \tilde u_\lambda$ by \ref{demo1}, because $\varphi$ has compact support.
    \item{(vi)} Now for the right side, $w_\rho$ strongly converges as $u_\rho$ does 
and  $f_\rho$ weakly converges to $f$.

To summarise \ref{demo5} converges to (remember that $\tilde P_\lambda^i$ are the coordinates of $\tilde P_\lambda$
which is horizontal)
$$
\int_{B_\infty(1)} \bigl(\tilde P_\lambda^jx_k - q^{jk} \tilde u_\lambda\bigr)X_j^\infty\varphi~d\mu_\infty %
= \int_{B_\infty(1)} fx_k \varphi~d\mu_\infty -\lambda \int_{B_\infty(1)}\varphi u_\lambda x_k~d\mu_\infty \numeq \label{demo6}
$$
furthermore if we put into equation \ref{demo1.1}, $\phi_\infty=\varphi x_k$ it gives
$$
 \int_{B_\infty(1)}fx_k\varphi~d\mu_\infty-\lambda\int_{B_\infty(1)}\varphi \tilde u_\lambda x_k~d\mu_\infty = \int_{B_\infty(1)}\tilde P_\lambda^j X_j^\infty(\varphi x_k)~d\mu_\infty  \numeq \label{equa}
 $$
and by mixing \ref{demo6} and \ref{equa} we get
for every $\varphi \in {\cal C}^\infty_c(B_\infty(1))$ the following equality : 
$$
\int_{B_\infty(1)}  \bigl(\tilde P_\lambda^j x_k - q^{jk} \tilde u_\lambda\bigr)X_j^\infty\varphi~d\mu_\infty =\int_{B_\infty(1)}\tilde P_\lambda^j X_j^\infty(\varphi x_k)~d\mu_\infty
$$
which in terms of distribution can be translated into :
$$
-\sum_{j=1}^{d_1}X_j^\infty\bigl(\tilde P_\lambda^j x_k - q^{jk} \tilde u_\lambda\bigr) =-\sum_{j=1}^{d_1} X_j^\infty\tilde P_\lambda^j~x_k
\iff  \tilde P_\lambda^k=\sum_{j=1}^{d_1}q^{jk}~X_j^\infty\tilde u_\lambda
$$
which allow us to conclude that $\tilde u_\lambda=G_\lambda f$. \qed

\paragraphe It is now easy to finish the proof, for theorem \ref{cvres} gives the compact
convergence of the resolvents fact which thanks to theorem \ref{cvssets} gives
the convergence of the spectrum which is what theorem \ref{theo1} is all about.

\partie{Macroscopic spectral rigidity}\label{example}
In this section we give the proof of the upper bound on the first eigenvalue and
of theorem \ref{thealst}.

\subpartie{The upper bound}
 
By proposition \ref{alinst} we have $B_\infty(1)\supset B_2(1)$ thus by the minmax property, 
we have
$$
\lambda_i(B_\infty(1)) \leq \lambda_i(B_2(1)) \numeq \label{egaliteDeslambda}
$$
for any $i$ and equality holds, by the maximum principle (see J.-M.~Bony \cite{bony}) if and only if
the two balls coincides.

\subpartie{Nilmanifolds having a one dimensional center}

The aim of this part is to characterize the case of equality in \ref{egaliteDeslambda}
for a class of nilmanifold, which contains the Heisenberg groups.

Let us first introduce the description of the metrics involved.

\proclaim Definition \theo. Let $N^{n+1}$ be a $2$-step nilmanifold whose kernel
is one dimensional. Suppose that there is a submersion $p$ of $N^{n+1}$
onto a flat torus $\T^n$. Let $(\alpha_1,\dots,\alpha_n)$ be the lift of an orthonormal base of harmonic
$1$-forms over the torus. 
Choose a $1$-form $\vartheta$ of $N^{n+1}$ such that $d\vartheta=p^*b$ where $b$
is a closed $2$-form over the torus (in other words we chose a connection). 
Let $g_\vartheta$ be the Riemanniann metric
such that the dual base of $(\alpha_1,\dots,\alpha_n,\vartheta)$ is orthonormal. Thus
$p$ becomes a Riemannian submersion. We will call such a metric {\rm pseudo left invariant}.%
\label{defpli}

The idea is that if your choice of connections gives a left invariant base of
vector fields then the  above construction gives a left invariant metric. Thus
this pseudo left invariant metric can be seen as perturbation of a left invariant metric,
by perturbation of a left invariant base of vector fields.

We are now able to give our precise claim

\proclaim Lemma \theo.
Let $(\H_{2n+1},g)$ be the $2n+1$-dimensional Heisenberg group, then its stable
norm coincides with its albanese Metric if and only if $g$ is pseudo left invariant.

which thanks to theorem \ref{theoup} allows us to obtain the following corollary :

\proclaim Corollary \theo.
Let $(\H_3,g)$ be the $3$-dimensional Heisenberg group, $B_g(\rho)$ the induced Riemannian
ball of radius $\rho$
on its universal cover and $\lambda_1\bigl(B_g(\rho)\bigr)$ the first
eigenvalue of the Laplacian on $B_g(\rho)$ for the Dirichlet problem.
\endgraf
Then
\item{1. } $\lim_{\rho\to +\infty} \rho^2\lambda_1\bigl(B_g(\rho)\bigr) = \lambda_1^\infty \leq \lambda_1^{{\cal H}_3}$
\item{2. } in case of equality the metric is pseudo left invariant.
\endgraf
Where $\lambda_1^{{\cal H}_3}$ is the first eigenvalue of the Kohn laplacian arising from
the Albanese metric on the unit ball of the carnot-caratheodory ball arising
from the same metric. \label{theoheisone}

In fact we are going to show a little more, we are going to focus on $2$-step nilmanifolds
for which the  center of their Lie algebra is of dimension $1$.

\proclaim Lemma \theo.
Let $(M^{n+1},g)$ be a $2$-step nilmanifold whose center is of dimension $1$,
then its stable norm and its Albanese metric coincides if and only if the metric
is pseudo left invariant.

\proof
The main idea comes from the fact  that
the albanese metric and the stable norm coincides if and only if every harmonic
$1$-form has constant length. 

Now take an orthonormal base of Harmonic $1$-forms $\alpha_1,\dots,\alpha_n$ (which can be seen
as lift of harmonic forms over the Albanese torus),
and consider their dual vector fields with respect to the metric $X_1,\dots,X_n$, they span ${\cal H}$.
Using the fact that for a closed one form $\alpha$ and any vector fields $X$ and $Y$ we have:
$$
\alpha([X,Y])=X\cdot\alpha(Y)-Y\cdot\alpha(X)
$$
we see that for any $i,j$ 
$$
[X_i,X_j] \in {\cal H}^\perp 
$$
(remark that we can also deduce from that fact that $[X_i,X_j]=2\nabla_{X_i}X_j$)
Now look at $Z$ the dual vector field to the $1$-form $Z^\flat=*(\alpha_1\land\dots\land\alpha_n)$ ($*$ is the Hodge
operator, thus this form is co-closed), its length is
constant by construction. Furthermore $Z$ belongs and spans ${\cal H}^\perp$.

For a co-closed one form $\alpha$ we have
$$
\sum_i \nabla_{X_i}\alpha\cdot\alpha_i + \nabla_Z\alpha\cdot Z^\flat =0
$$
which implies that for $i=1,\dots, n$
$$
[X_i,Z] = 0
$$
which also implies that $Z$ is a killing field.

hence we have
$$\eqalign{
[X_i,X_j] &= f_{ij} Z }
$$
where $f_{ij}$ are some functions, which are not all zero, otherwise
our nilmanifold will be locally a flat torus.

Let us remark that ($Z$ being a killing field) we have
$$
dZ^\flat(X,Y)=2g(\nabla_XZ,Y) \numeq \label{crochet4}
$$
thus if we decompose $dZ^\flat$ in the base given by $\alpha_i\land\alpha_j$, $Z^\flat\land\alpha_i$ and $Z^\flat\land Z^\flat$ for all $i,j$ then
thanks to \ref{crochet3} and \ref{crochet4} we get that
$$
dZ^\flat=\sum_{i<j} f_{ij} \alpha_i\land\alpha_j
$$
In other words $dZ^\flat$ is horizontal and 
there exists some $2$-form $\beta$ on the albanese torus such that
$dZ^\flat=\pi^*\beta$.

Hence $g$ is pseudo left invariant.


\qed

\subpartie{On the higher dimensional case}
The case of tori, which is quite exceptional has been treated in \cite{ver1}.

For other, non abelian nilmanifolds,
it seems that things are not that simple, as in the case of the Heisenberg groups,
it seems that there can be
metrics which are not left invariant, but for which the equality between the
stable norm and the albanese metric holds.

The interested reader will be able to consult \cite{pac}, where this question is
studied in details, and where we explain the rigidity involved.

\subpartie{Graded nilmanifolds  with totally geodesic fibers \\ over a Torus}

There is one last particular case we would like to study, the case where the nilmanifold
is graded (i.e. its algebra is nilpotent and graded as defined in section \ref{graded}) and
the metric on $(M^n,g)$ is as follows : We suppose that the first betti number $b_1(M^n)=k$,
we recall that ${\cal H}$ is the horizontal distribution coming from $V_1$ 
(see sections \ref{horizontald} and \ref{graded}). Moreover we assume that
we have the following riemanniann submersion, with totally geodesics fibers
and with a metric equivariant on the fibers :
$$
[M,M] \hookrightarrow (M^n,g)\ {\buildrel p\over \longrightarrow}\ (\T^k,\check g)
$$
where $dp_x$ is an isometry (we write $\hat g=g_{\mid {\cal H}}$) 
from $({\cal H}_x,\hat g_x)$ to $(T_{p(x)}\T^k,\check g_{p(x)})$.

Then, in case of equality in the theorems \ref{theoup} and \ref{theolast}, the albanese
map is a riemanniann submersion, which implies that $\check g$ is flat. Which in turn, using
our assumptions implies that the metric $g$ is left invariant (indeed see chapter 9 section F in 
\cite{besse}).  In other words :

\proclaim Proposition \theo.
Let $(M,g)$ satisfying the above assumptions. 
The albanese metric and the stable norm coincides if and
only if the metric is left invariant.

In other words, we could say heuristically that for sub-riemannian metrics the equality case
in theorems \ref{theoup} (which holds in that context too, see \cite{vernicos}
for the convergence of the spectrum)
characterises the left-invariant sub-riemannian metrics.


\bibliography


\bibitem[Ale02]{alexopoulos}
\bgroup\bf G.~K. Alexopoulos\egroup{}.
\newblock {\em Sub-Laplacians with Drift on Lie Groups of Polynomial Volume
  Growth}, volume 155, number 739 of {\em Memoirs of the AMS}.
\newblock AMS, January 2002.
\finitem

\bibitem[BI95]{bi2}
\bgroup\bf D.~Burago\egroup{} and \bgroup\bf S.~Ivanov\egroup{}.
\newblock On asymptotic volume of tori.
\newblock {\em GAFA}, 5(5):800--808, 1995.
\finitem

\bibitem[BLP78]{blp}
\bgroup\bf A.~Bensousan\egroup{}, \bgroup\bf J.-L. Lions\egroup{}, and
  \bgroup\bf G.~Papanicolaou\egroup{}.
\newblock {\em Asymptotic analysis for periodic structures}.
\newblock Studies in mathematics and its applications. North Holland, 1978.
\finitem

\bibitem[BMT96]{bmt1}
\bgroup\bf M.~Biroli\egroup{}, \bgroup\bf U.~Mosco\egroup{}, and \bgroup\bf
  N.~A. Tchou\egroup{}.
\newblock Homogenization for degenerate operators with periodical coefficients
  with respect to the heisenberg group.
\newblock {\em C. R. Acad. Sci. Paris}, t. 322, S\'erie I:439--444, 1996.
\finitem

\bibitem[BMT97]{bmt2}
\bgroup\bf M.~Biroli\egroup{}, \bgroup\bf U.~Mosco\egroup{}, and \bgroup\bf
  N.~A. Tchou\egroup{}.
\newblock Homogenization by the heisenberg group.
\newblock {\em Advances in Mathematics}, 7:809--837, 1997.
\finitem

\bibitem[Bon69]{bony}
\bgroup\bf J.-M.~Bony\egroup{}.
\newblock Principe du maximum, in\'egalite de {H}arnack et unicit\'e du
  probl\`eme de {C}auchy pour les op\'erateurs elliptiques d\'eg\'en\'er\'es.
\newblock {\em Ann. Inst. Fourier (Grenoble)}, 19(fasc. 1):277--304 xii, 1969.
\finitem

\bibitem[BP92]{bepe}
\bgroup\bf R.~Benedetti\egroup{} and \bgroup\bf C.~Petronio\egroup{}.
\newblock {\em Lectures on Hyperbolic Geometry}.
\newblock Springer-Verlag, 1992.
\finitem

\bibitem[Bes87]{besse}
\bgroup\bf A.~L. Besse\egroup{}.
\newblock {\em Einstein Manifolds}.
\newblock Springer-Verlag, 1987.
\finitem

\bibitem[Bro85a]{brooks}
\bgroup\bf R.~Brooks\egroup{}.
\newblock The bottom of the spectrum of a {R}iemannian covering.
\newblock {\em J. Reine Angew. Math.}, 357:101--114, 1985.
\finitem


\bibitem[Fed69]{federer}
\bgroup\bf H.~Federer\egroup{}.
\newblock {\em Geometric Measure Theory}.
\newblock Springer Verlag, 1969.
\finitem

\bibitem[GHL90]{ghl}
\bgroup\bf S.~Gallot\egroup{}, \bgroup\bf D.~Hulin\egroup{}, and \bgroup\bf
  J.~Lafontaine\egroup{}.
\newblock {\em Riemannian geometry}.
\newblock Universitext. Springer-Verlag, second edition, 1990.
\finitem

\bibitem[Gro81]{gromov}
\bgroup\bf M.~Gromov\egroup{}.
\newblock Groups of polynomial growth and expanding maps.
\newblock {\em Inst. Hautes \'etudes Sci. Publ. Math.}, (53):53--73, 1981.
\finitem


\bibitem[KS]{kushi}
\bgroup\bf K.~Kuwae\egroup{} and \bgroup\bf T.~Shioya\egroup{}.
\newblock Convergence of spectral structures: a functional analytic theory and
  its applications to spectral geometry.
\newblock preprint.
\finitem


\bibitem[Mas93]{Dmaso}
\bgroup\bf Dal Maso\egroup{}.
\newblock {\em An Introduction to $\Gamma$-convergence.}
\newblock Birkh{\"a}user, 1993.
\finitem


\bibitem[Mos94]{mosco}
\bgroup\bf U.~Mosco\egroup{}.
\newblock Composite media and asymptotic dirichlet forms.
\newblock {\em J. Funct. Anal.}, 123(2):368--421, 1994.
\finitem


\bibitem[NV]{pac}
\bgroup\bf P.A.~Nagy\egroup{} and \bgroup\bf C.~Vernicos\egroup{}.
\newblock On manifolds with eigenforms of constant length.
\newblock In preparation.
\finitem

\bibitem[Nom54]{nomizu}
\bgroup\bf K.~Nomizu\egroup{}.
\newblock On the cohomology of compact homogeneous spaces of nilpotent lie
  groups.
\newblock {\em Annals of Math.}, 59(3):531--538, 1954.
\finitem

\bibitem[Pan82]{pansu}
\bgroup\bf P.~Pansu\egroup{}.
\newblock {\em G\'eometrie du groupe de Heisenberg}.
\newblock Th\`ese de docteur 3\`eme cycle, Universit\'e Paris VII, 1982.
\finitem

\bibitem[Pan99]{pansu2}
\bgroup\bf P.~Pansu\egroup{}.
\newblock Profil isop\'erim\'etrique, m\'etriques p\'eriodiques et
  formes d'\'equilibre des cristaux.
\newblock pr\'epublication d'orsay, 1999.
\finitem

\bibitem[Rud91]{rudinfa}
\bgroup\bf W.~Rudin\egroup{}.
\newblock {\em Functional Analysis}.
\newblock International series in Pure and Applied Mathematics. Mc Graw-Hill,
  second edition, 1991.
\finitem

\bibitem[Sun89]{sunada}
\bgroup\bf T.~Sunada\egroup{}.
\newblock Unitary representations of fundamental groups and the spectrum of
  twisted {L}aplacians.
\newblock {\em Topology}, 28(2):125--132, 1989.
\finitem

\bibitem[Ver01]{vernicos}
\bgroup\bf C.~Vernicos\egroup{}.
\newblock {\em Spectres asymptotiques des nilvari\'et\'es gradu\'ees}.
\newblock Th\`ese de doctorat, Universit\'e Grenoble I, Joseph Fourier,
  2001.
\finitem

\bibitem[Ver02]{ver1}
\bgroup\bf C.~Vernicos\egroup{}.
\newblock The macroscopical sound of tori.
\newblock preprint.
\finitem

\bibitem[ZKON79]{gconv}
\bgroup\bf V.V. Zhikov\egroup{}, \bgroup\bf S.M. Kozlov\egroup{}, \bgroup\bf
  O.A. Oleinik\egroup{}, and \bgroup\bf Kha~T'en Ngoan\egroup{}.
\newblock Averaging and g-convergence of differential operators.
\newblock {\em Russian Math. Surveys}, 34(5):69--147, 1979.
\finitem
\bibliend

\bigskip

\vbox{\hsize 7cm  \parindent=0cm \parskip=0cm
{\bf Constantin Vernicos}\par
\sf Universit\'e de Neuch\^atel\par
Institut de Math\'ematiques\par
11, rue Emile Argand\par
2007 Neuch\^atel\par
Switzerland\par
mail: \tt Constantin.Vernicos@unine.ch}

\bye